# A REMARK ON MAPPING TORI OF FREE GROUP ENDOMORPHISMS

ILYA KAPOVICH

If $\phi : G \to G$ is an endomorphism (not necessarily injective) of a group $G$ then the *mapping torus group* $M(\phi)$ is defined as:

$$M(\phi) = \langle G, t | t^{-1}gt = \phi(g), g \in G \rangle.$$

We observe the following:

**Proposition 0.1.** *Let $F$ be a finitely generated free group and let $\phi : F \to F$ be an endomorphism. Then there exists a finitely generated free group $F_1$ and an injective endomorphism $\psi : F_1 \to F_1$ such that $M(\phi) \cong M(\psi)$.*

*Proof.* It is well-known that the kernels of the powers of $\phi$ stabilize (see for example [3]), that is, there exists $k > 0$ such that $ker(\phi^k) = ker(\phi^n)$ for all $n \geq k$. (This easily follows from the stabilization of ranks of the free groups $\phi^k(F)$ and from Hopficity of finitely generated free groups.) Put $N = ker(\phi^k)$. Then $\phi$ factors through to an injective endomorphism $\overline{\phi} : F/N \to F/N$. The group $F/N$ is isomorphic to $\phi^k(F) \leq F$, so that $F/N$ is a free group of finite rank.

Let

$$G = M(\phi) = \langle F, t | tft^{-1} = \phi(f), f \in F \rangle.$$

If $f \in N$ then $t^k f t^{-k} = \phi^k(f) = 1$ and hence $f = 1$ in $G$. We perform a Tietze transformation on the presentation of $G$ above and add the relations $f = 1$ for all $f \in N$. This produces a presentation of $G$ as the HNN-extension $M(\overline{\phi})$ of the free group $F/N$ along an injective endomorphism $\overline{\phi}$. □

M. Feign and M. Handel [1] proved that mapping tori of injective endomorphisms of finitely generated free groups are coherent (that is all their finitely generated subgroups are finitely presentable). Later R. Geoghegan, M. Mihalik, M. Sapir and D. Wise [2] used this result to establish that such mapping tori are Hopfian. Together with Proposition 0.1 this implies that the requirement of injectivity can be dropped in both cases:

**Corollary 0.2.** *Let $\phi : F \to F$ be an endomorphism of a finitely generated free group $F$. Then $M(\phi)$ is coherent and Hopfian.*

---

Department of Mathematics, University of Illinois at Urbana-Champaign, 1409 West Green Street, Urbana, IL 61801, USA

*E-mail address*: kapovich@math.uiuc.edu